\documentclass[12pt,twoside]{article}
\usepackage{amsmath,amssymb}

\title{On the existence and multiplicity of positive solutions to classes of steady state reaction diffusion systems with multiple parameters} %
\author{A. Shabanpour $^{a}$, S.H. Rasouli $^{a}$, N. Fonseka $^{b}$ \\
$^{a}$ Department of Mathematics, Faculty of Basic Sciences,\\
Babol Noshirvani University of Technology, Babol, Iran\\e-mail:
javadshabanpour50@gmail.com (A. Shabanpour)\\ s.h.rasouli@nit.ac.ir (S.H. Rasouli)\\
$^{b}$ School of Arts and Sciences, Carolina University,\\ Winston-Salem, NC 27101, USA\\
e-mail: fonsekan@carolinau.edu }
\date{}
\begin{document}
\maketitle
\begin{center}
{\bf\large Abstract}\\
\end{center}
we study positive solutions to the steady state reaction diffusion systems of the form:
\begin{equation}
\left\{\begin{array}{ll}
-\Delta u  = \lambda f(v)+\mu h(u), &   \Omega,\\
-\Delta v  =  \lambda g(u)+\mu q(v),& \Omega,\\
\frac{\partial u}{\partial \eta}+\sqrt[]{\lambda +\mu}\, u=0,& \partial\Omega,\\
\frac{\partial v}{\partial \eta}+\sqrt[]{\lambda +\mu}\, v=0, & \partial\Omega,\\
\end{array}\right.
\end{equation}
where ${\lambda,\mu>0}$ are positive parameters, ${\Omega}$ is a bounded in  $\mathbb{R}^{N}$$(N>1)$ with smooth boundary ${\partial \Omega}$, or ${\Omega=(0,1)}$, ${ \frac{\partial z}{\partial \eta} }$ is the outward normal derivative of $z$. Here $f, g, h, q\in C^{2} [0,r)\cap C[0,\infty)$ for some $r>0$. Further, we assume that $f, g, h$ and $q$ are increasing functions such that $f(0) = g(0) = h(0) = {q}(0) = 0$, $f^\prime(0), g^\prime(0), h^\prime(0), q^\prime(0) > 0$, and  $\lim\limits_{s\to \infty}\frac{f(M g(s))}{s}=0$ for all $M>0$. Under certain additional assumptions on $f, g, h$ and $ q$ we prove our existence and multiplicity results. Our existence and multiplicity results are proved using sub-super solution methods.\\\\

\hspace{-0.6 cm}Keywords:  Positive solutions; Reaction diffusion systems; Sub-super solutions.
\section{Introduction}
\hspace{0.6 cm}Many papers have been devoted during recent years to the study of
reaction-diffusion systems which arise very often in applications in fields such as, mathematical
biology, chemical reactions, and combustion theory (see \cite{aa-nf-rs1}, \cite{aa-nf-rs}, \cite{ja-mr-rs}, and \cite{ramaswami}).\\

In this paper, we consider the following reaction-diffusion system with multiple parameters
\begin{equation}
\left\{\begin{array}{ll}
-\Delta u  = \lambda f(v)+\mu h(u), &   \Omega,\\
-\Delta v  =  \lambda g(u)+\mu q(v),& \Omega,\\
\frac{\partial u}{\partial \eta}+\sqrt[]{\lambda +\mu}\, u=0,& \partial\Omega,\\
\frac{\partial v}{\partial \eta}+\sqrt[]{\lambda +\mu}\, v=0, & \partial\Omega,\\
\end{array}\right.
\end{equation}
where ${\lambda,\mu>0}$ are positive parameters, ${\Omega}$ is a bounded in  $\mathbb{R}^{N}$$(N>1)$ with smooth boundary ${\partial \Omega}$, or ${\Omega=(0,1)}$, ${ \frac{\partial z}{\partial \eta} }$ is the outward normal derivative of $z,$ and $f, g, h, q$ are continuous increasing functions.\\

A lot of phenomena
appearing in ecology are accompanied with diffusion and most of such phenomena can be described in the form of partial differential equations with diffusion
terms. Among them, reaction-diffusion equations and their related systems have
attracted interests of many researchers and have been investigated very extensively.
Moreover, new mathematical models with nonlinear diffusion have been proposed
in recent years; so that the importance of the study of reaction-diffusion equations
are increasing from the ecological view-point as well as the mathematical one.\\

On the other hand, in mathematical modeling, elliptic partial differential equations are used together with
boundary conditions specifying the solution on the boundary of the domain. Commonly, problems with boundary conditions that are linear functions of the values or
normal derivatives of the solutions on the boundary have been studied extensively
(see \cite{rsc-cc,rc-hl-lm-js,pd-shr,ssrf-htt,jz} and references therein).\\

Our work is motivated by the results in the literature for single-equation case, namely equations of the form
\begin{equation}
\left\{\begin{array}{ll}
-\Delta u  = \lambda f(u);  &   \Omega,\\
\frac{\partial u}{\partial \eta}+\mu (\lambda) u=0; & \partial\Omega,\\
\end{array}\right.
\end{equation}
where $f\in C^{2}[0,\infty),$ and $\mu\in C[0,\infty)$ is strictly such that $\mu(0)\geq 0.$ In recent history, there has been growing interest on the study of the reaction diffusion models where the parameter influences the equation as well as the boundary condition (see \cite{aa-nf-jq-rs,aa-nf-rs1,aa-nf-rs,jtc-nf-jg-jl-rs,nf-jg-rs-bs,nf-jm-rs,nf-am-rs-bs,nf-rs-bs-ks,nf-rs-jg-qm-bs}). See \cite{nf-rs-bs-ks}, where the authors discussed the existence, nonexistence, multiplicity, uniqueness results to the problem $(3).$  Recently, in \cite{aa-nf-rs1,aa-nf-rs} the authors extended the results in \cite{nf-rs-bs-ks} to the following system (see \cite{aa-nf-rs} for detail)
\begin{equation}
\left\{\begin{array}{ll}
-\Delta u  = \lambda f(v) ;  &   \Omega,\\
-\Delta v  =  \lambda g(u) ;& \Omega,\\
\frac{\partial u}{\partial \eta}+\sqrt[]{\lambda }\, u=0; & \partial\Omega,\\
\frac{\partial v}{\partial \eta}+\sqrt[]{\lambda } \,v=0; & \partial\Omega,\\
\end{array}\right.
\end{equation}
where $\lambda >0$ is positive parameter, $f$ and $g$ are continuous functions such that $f(0) = 0 = g(0)$ and $\lim \limits_{x\to \infty}\frac{f(M g(x))}{x} = 0$ for all $M>0$. Here we focus on further extending the study in \cite{aa-nf-rs} to the system $(1)$ which features multiple parameters and stronger coupling.

\section{Our results}
We first introduce our hypotheses $(H_{1})-(H_{3})$.
\begin {description}
\item{\bf (H1)} $f, g, h, q \in C[0,\infty)$  increasing functions such that $f(0) = g(0) = h(0) = q(0)=0;$
\end {description}
\begin {description}
\item{\bf (H2)} $\lim\limits_{x\to \infty}\frac{f(M g(x))}{x}=0$ for all $M>0$ (combined sublinear effect at infinity);
\end {description}
\begin {description}
\item{\bf (H3)} There exists $a>0$ such that $f, g, h, q \in C^{1}[0 ,a)$ and $f^\prime(0),g^\prime(0),h^\prime(0), q^\prime(0)>0$.
\end {description}

\noindent Next we introduce the hypotheses $(F_{1})-(F_{3})$.

\begin {description}
\item{\bf (F1)} $f, g, h, q \in C[0,\infty)$  are bounded functions
\end {description}

\begin {description}
\item{\bf (F2)} $h$ is a function that is sublinear at infinity (ie. $\lim\limits_{s\rightarrow\infty}\frac{h(s)}{s}=0$),  there exist $\gamma, \beta >0$ such that  $\gamma \beta<1$, $g(s)\leq s^\gamma$ for $s\gg 1$, $q(s)\leq s^\beta$ for $s\gg 1$, and $\lim\limits_{s\rightarrow\infty}g'(s)>0$
\end {description}

\begin {description}
\item{\bf (F3)} $q$ is a function that is sublinear at infinity (ie. $\lim\limits_{s\rightarrow\infty}\frac{q(s)}{s}=0$),  there exist $\gamma, \beta >0$ such that  $\gamma \beta<1$, $f(s)\leq s^\gamma$ for $s\gg 1$, $h(s)\leq s^\beta$ for $s\gg 1$, and $\lim\limits_{s\rightarrow\infty}f'(s)>0$
\end {description}

Without loss of generality, we assume that  $f^\prime(0)\geq g^\prime(0)$. Further, we assume that one of $(F_{1}), (F_{2}), (F_{3})$ holds throughout the paper.\\

 Here we introduce an eigenvalue $A_{1}$, which we use to prove our theorems. To accomplish this goal, recall recent results for the eigenvalue problem
\begin{equation}
\left\{\begin{array}{ll}
-\Delta \Theta = K  \Theta;  &   \Omega,\\
\frac{\partial \Theta}{\partial \eta}+ \tau \sqrt{K}\Theta =0; & \partial\Omega.   \\
\end{array}\right.
\end{equation}
where $\tau > 0$. Namely, let $K_{1}(\tau ) $ be its principal eigenvalue (see \cite{jg-qm-sr-rs}), and let $\Theta$ be the corresponding normalized positive eigenfunction of (5). Now consider the eigenvalue problem
\begin{equation}
\left\{\begin{array}{ll}
-\Delta \Gamma = \overline{K} g^\prime(0) \Gamma;  &   \Omega,\\
\frac{\partial \Gamma}{\partial \eta}+  \sqrt{\overline{K}}\Gamma =0; & \partial\Omega.  \\
\end{array}\right.
\end{equation}

Noting that the substitution $K = \overline{K} g^\prime(0) $  reduces (6) to (5), we easily see that the principal
eigenvalue of (6) is $\frac{K_{1}(\tau )  }{ g^\prime(0) }$ with $\tau= \frac{1}{\sqrt{ g^\prime(0)} }$. Define
\begin{center}
$A_{1}:=\frac{K_{1}(\tau )  }{ g^\prime(0) }$
\end{center}

Next, for a given $\tilde{\lambda}$, we recall the for the eigenvalue problem
\begin{equation}
\left\{\begin{array}{ll}
-\Delta \Theta = ( \varrho + \tilde{\lambda}) \Theta;  &   \Omega,\\
\frac{\partial \Theta}{\partial \eta}+  \tau \sqrt{\tilde{\lambda}} \Theta=0; & \partial\Omega.   \\
\end{array}\right.
\end{equation}

Denoting by $\varrho_{\tilde{\lambda}}  (\tilde{\lambda}, \tau) $ its principal eigenvalue and by $\phi_{\tilde{\lambda}}>0$ the corresponding eigenfunction of $(7)$ such that ${\| \phi_{\tilde{\lambda}}\|}_\infty=1$, the following results were established in \cite{jg-qm-sr-rs}:

\begin{equation}
\left\{\begin{array}{ll}
 \varrho_{\tilde{\lambda}} > 0;  & \tilde{\lambda}< K_{1}(\tau),\\
 \varrho_{\tilde{\lambda}}= 0; & \tilde{\lambda}= K_{1}(\tau) ,  \\
 \varrho_{\tilde{\lambda}}< 0;  & \tilde{\lambda}> K_{1}(\tau). \\
\end{array}\right.
\end{equation}

Hence, by the substitution $\lambda g^\prime(0)=  \tilde{\lambda}$, denoting by $ \varrho_{\lambda}$ the principal eigenvalue of
\begin{equation}
\left\{\begin{array}{ll}
-\Delta \Gamma = (\varrho+\lambda g^\prime(0)) \Gamma;  &   \Omega,\\
\frac{\partial \Gamma}{\partial \eta}+  \sqrt{\lambda}\Gamma =0; & \partial\Omega, \\
\end{array}\right.
\end{equation}
then, the following results were established in \cite{aa-nf-rs}:
\begin{equation}
\left\{\begin{array}{ll}
 \varrho_{\lambda}> 0;  & {\lambda}< A_{1},\\
 \varrho_{\lambda}= 0; & {\lambda}= A_{1},  \\
 \varrho_{\lambda}< 0;  & {\lambda}> A_{1}. \\
\end{array}\right.
\end{equation}

To precisely state our existence results we consider the unique classical positive solution $\xi$ of the following  elliptic problem
\begin{equation}
\left\{\begin{array}{ll}
-\Delta \xi =1;  &   \Omega,\\
\frac{\partial \xi}{\partial \eta}+ \xi =0; & \partial\Omega.\\
\end{array}\right.
\end{equation}

Further, for positive constants $a,b$ such that $a<b,$ define
$$
Q_{1}(a):=\min\{\frac{a}{f(a)}, \frac{a}{g(a)},\frac{a}{h(a)},\frac{a}{q(a)}\},
$$
and
$$
Q_{2}(b):=\max \{\frac{b}{f(b)}, \frac{b}{g(b)},\frac{b}{h(b)},\frac{b}{q(b)}\}.\\
$$

Now we state our  results.\\\\
{\bf Theorem 2.1.} Assume $(H_{1})-(H_{3})$  hold. Let $R$ be the radius of the largest inscribed ball in $\Omega$, then for $0<\epsilon<R,$ we define
$$
C_1=\inf_{\epsilon} \frac {NR^N-1 }{   \epsilon^N R-\epsilon}.
$$
Then the system $(1)$ has a positive solution ($u_{\lambda , \mu}$,$v_{\lambda , \mu}$) for $\lambda + \mu >A_{1}$, such that ${\|u_{\lambda , \mu}\|}_\infty,$ ${\|v_{\lambda , \mu}\|}_\infty$ $\rightarrow\infty$ as ${\lambda + \mu}$ $\rightarrow\infty$. Further, if $\frac{Q_{1}(a)}{Q_{2}(b)}>2C_1{\|\xi\|}_\infty$ and  $Q_{1}(a)>\max 2\{A_1,1\}{\|\xi\|}_\infty,$ then system $(1)$ has at least three positive solutions for $\lambda ,\mu \in (\max \{A_1, C_{1}(\Omega)Q_{2}(b),1\},\frac{Q_{1}(a)}{{2\|\xi\|}_\infty}).$\\\\
{\bf Theorem 2.2.} If there exists $r>0$ such that $f, g \in C^{2} [0,r)$, $f(0) = g(0) = h(0) = {q}(0) = 0$, $f^\prime(0) = g^\prime(0) = h^\prime(0) = q^\prime(0) > 0$, and $ f^{\prime\prime}(s) , g^{\prime\prime}(s), h^{\prime\prime}(s) , q^{\prime\prime}(s) < 0$ for all $ s\in [0, r)$, then system $(1)$ has a positive solution $(u_{\lambda,\mu} , v_{\lambda,\mu})$ for $\lambda+\mu>A_1 $ and $\lambda +\mu \approx A_1$ such that  ${\|  u _{\lambda ,  \mu}  \|}$, ${\|  v _{\lambda ,  \mu}  \|} \rightarrow 0 $ as ${\lambda +  \mu}  \rightarrow A_1^+$.\\

The reminder of this paper is organized as follows. In Section $3,$ we present some preliminaries. In Section $4,$  using the sub-supersolutions technique, we prove Theorems $2.1$ and $2.2.$ In section 5, we discuss an example.

\section{Preliminaries and the method of sub-supersolutions}
\hspace{0.6 cm}In this section, we state method sub-supersolution (see \cite{gaa-shr}). A pair of nonnegative
functions $(\psi_{1},\psi_{2})\in \Big(W^{1,2}\cap
C(\overline{\Omega})\Big)\times  \Big(W^{1,2}\cap C(\overline{\Omega})\Big)$ and
$(z_{1},z_{2})\in \Big(W^{1,2}\cap C(\overline{\Omega})\Big)\times \Big(W^{1,2}\cap
C(\overline{\Omega})\Big)$ are called a subsolution and supersolution of
$(1)$ if they satisfy
\begin{equation}
\left\{\begin{array}{ll}
-\Delta \psi_{1} \leq\lambda f(\psi_{2})+\mu h(\psi_{1}) ;  &   \Omega,\\
-\Delta \psi_{2} \leq \lambda g(\psi_{1})+\mu q(\psi_{2}) ;& \Omega,\\
\frac{\partial  \psi_{1}}{\partial \eta}+\sqrt[]{\lambda +\mu}  \psi_{1}\leq0; & \partial\Omega,\\
\frac{\partial  \psi_{2}}{\partial \eta}+\sqrt[]{\lambda +\mu}  \psi_{2}\leq0; & \partial\Omega,\\
\end{array}\right.
\end{equation}
and
\begin{equation}
\left\{\begin{array}{ll}
-\Delta z_{1} \geq\lambda f(z_{2})+\mu h(z_{1}) ;  &   \Omega,\\
-\Delta z_{2} \geq \lambda g(z_{1})+\mu q(z_{2}) ;& \Omega,\\
\frac{\partial z_{1}}{\partial \eta}+\sqrt[]{\lambda +\mu} z_{1} \geq0; & \partial\Omega,\\
\frac{\partial z_{2}}{\partial \eta}+\sqrt[]{\lambda +\mu}  z_{2}\geq0; & \partial\Omega,\\
\end{array}\right.
\end{equation}
respectively. It is well known that if there exist sub and
supersolutions $(\psi_{1},\psi_{2})$ and $(z_{1},z_{2})$
respectively of $(1)$ such that $(\psi_{1},\psi_{2})\leq
(z_{1},z_{2}),$ then $(1)$ has a solution $(u,v)$ such that
$$
(u,v)\in[(\psi_{1},\psi_{2}),(z_{1},z_{2})] (see \cite {ha, fi, pao}).
$$
By strict sub and super-solutions mean functions $(\psi_{1},\psi_{2})$ and $(z_{1},z_{2})$ for which  inequalities in $(12)$ and $(13)$ are strict.\\

Our multiplicity results are obtained by constructing sub and super-solution pairs that satisfy the following Lemma:\\\\
{\bf Lemma 3.1.} (See \cite{fi,rs}). Suppose the system $(1)$
has a sub-solution $( \psi_{1},\psi_{2}),$ a strict super-solution
$(\zeta_{1},\zeta_{2}),$ a strict sub-solution $(w_{1},w_{2}),$ and
a super-solution $(z_{1},z_{2})$ for $(1)$ such that
$$
( \psi_{1},\psi_{2})\leq (\zeta_{1},\zeta_{2})\leq (z_{1},z_{2}),
$$
$$
( \psi_{1},\psi_{2})\leq (w_{1},w_{2})\leq (z_{1},z_{2}),
$$
and $(w_{1},w_{2})\nleq (\zeta_{1},\zeta_{2}).$ Then $(1)$ has at least three distinct positive solutions $(u_{i},v_{i}),$ $i = 1, 2, 3$ such
that
$$
 (u_{1}, v_{1})\in [( \psi_{1},\psi_{2}), (\zeta_{1},\zeta_{2})],\,\,\, (u_{2}, v_{2})\in [( w_{1},w_{2}), (z_{1},z_{2})]
$$
and
$$
(u_{3}, v_{3}) \in \Big[( \psi_{1},\psi_{2}),
(z_{1},z_{2})\Big]\backslash \Big(\Big[( \psi_{1},\psi_{2}),
(\zeta_{1},\zeta_{2})\Big]\cup \Big[ ( w_{1},w_{2}),
(z_{1},z_{2})\Big] \Big).
$$

\section{Proofs of our existence and multiplicity results}

In this section, we prove our existence and multiplicity results.

{\bf Proof of Theorem 2.1.} First we show the existence of positive solution for $\lambda + \mu>A_{1}$. Let  ($\psi_{1} $,$\psi_{2} $) = ($m\phi_{\lambda ,\mu } $, $m\phi_{\lambda ,\mu  } $), where $\phi_{\lambda , \mu} $ is the normalized positive eigenfunction of

\begin{equation}
\left\{\begin{array}{ll}
-\Delta \phi =(\varrho_{\lambda,\mu}+ (\lambda +\mu)  g^\prime(0))\phi ;  &   \Omega,\\
\frac{\partial \Gamma}{\partial \eta}+\sqrt[]{\lambda +\mu}  \phi=0; & \partial\Omega.\\
\end{array}\right.
\end{equation}
Define $H(s)=(\varrho_{\lambda ,\mu   }+ (\lambda +\mu)   g^\prime(0)) s -  (\lambda +\mu)   f(s)$. Then we have $H(0)= 0$ and $H^\prime(0)  \leq 0$
since $ \varrho_{\lambda ,\mu   }< 0$ for $\lambda +\mu>A_{1}$ and $f^\prime(0)\geq g^\prime(0)$. This implies that $H(s) \leq 0$ for $s \approx 0.$  Thus for  $m \approx 0$ we have
$$
(\varrho_{\lambda ,\mu   }+ (\lambda +\mu)   g^\prime(0))m\phi_{\lambda ,\mu} - \lambda f(m\phi_{\lambda ,\mu}) \leq 0.
$$
Note that $-\mu h( m\phi_{\lambda ,\mu}) \leq 0.$ Therefore, we have
$$
 (\varrho_{\lambda,\mu }+  (\lambda +\mu)   g^\prime(0))m\phi_{\lambda,\mu} - \lambda f(m\phi_{\lambda ,\mu}) - \mu h(m\phi_{\lambda,\mu }) \leq 0.
$$

By $(14)$ we have

\begin{center}
 $-\Delta \psi_{1} \leq\lambda f(\psi_{2})+\mu h(\psi_{1})$ for $m \approx 0.$
\end{center}
Similarly, by analyzing $\overline{H}(s)=(p+ (\lambda +\mu)   g^\prime(0)) s - (\lambda +\mu)q(s)$ one can show that
\begin{center}
 $-\Delta \psi_{2} \leq\lambda g(\psi_{1})+\mu q(\psi_{2}) $ for $m \approx 0.$
\end{center}

Further, on the boundary, we have

\begin{center}
$\frac{\partial  \psi_{1}}{\partial \eta}+\sqrt[]{\lambda +\mu}  \psi_{1} = 0$ and $ \frac{\partial  \psi_{2}}{\partial \eta}+\sqrt[]{\lambda +\mu}  \psi_{2} = 0$.
\end{center}

Thus  ($\psi_{1} $,$\psi_{2} $) is a subsolution of $(1)$ for  $m \approx 0.$\\

Now let $\zeta _{\lambda,\mu}$ be the unique positive solution of

\begin{equation}
\left\{\begin{array}{ll}
-\Delta \zeta =1;  &   \Omega,\\
\frac{\partial \zeta}{\partial \eta}+\sqrt[]{\lambda +\mu}\, \zeta =0; & \partial\Omega. \\
\end{array}\right.
\end{equation}

 We consider three different cases:

$(i)$ \noindent  Assume that (F1) is satisfied:
Let

$$
  (z_{1}, z_{2}) = \Big( (\lambda +\mu)   M_{ \lambda ,\mu   }\frac{\zeta_{\lambda,\mu }}{{\|\zeta_{\lambda,\mu}\|}_\infty},  (\lambda +\mu)    M_{\lambda,\mu }\frac{\zeta_{\lambda ,\mu}}{{\|\zeta_{\lambda ,\mu}\|}_\infty}\Big)
$$
Where $M_{\lambda,\mu}$ is a positive quantity. We choose $M_{\lambda,\mu }\gg1$ such that
\begin{center}
$\frac{M_{\lambda ,\mu} (\lambda +\mu)   }{{\|\zeta_{\lambda ,\mu}\|}_\infty}\geq \lambda $f\Big($\frac{  (\lambda +\mu)   M_{\lambda,\mu }\zeta_{\lambda ,\mu}}{{\|\zeta_{\lambda,\mu }\|}_\infty}$\Big) + $\mu h\Big(\frac{  (\lambda +\mu)   M_{\lambda ,\mu}\zeta_{\lambda ,\mu}}{{\|\zeta_{\lambda,\mu }\|}_\infty}\Big).$
\end{center}
Then we have
\begin{eqnarray*}
\begin{split}
-\Delta z_{1}&=  \frac{ (\lambda +\mu)   M_{\lambda,\mu }}{{\|\zeta_{\lambda ,\mu}\|}_\infty} \\
&\geq \lambda f\Big(\frac{ (\lambda +\mu) M_{\lambda,\mu }\zeta_{\lambda,\mu }}{{\|\zeta_{\lambda,\mu }\|}_\infty}\Big) + \mu h\Big(\frac{ (\lambda +\mu) M_{\lambda,\mu }\zeta_{\lambda,\mu }}{{\|\zeta_{\lambda,\mu }\|}_\infty}\Big).\\
\end{split}
\end{eqnarray*}

This implies

\begin{center}
$-\Delta z_{1} \geq{\lambda f(z_{2})+\mu h(z_{1})}$ for $M_{\lambda,\mu }\gg1$.
\end{center}
Similarly, $-\Delta z_{2} \geq \lambda g(z_{1})+\mu q(z_{2}) $ for $M_{\lambda,\mu }\gg1$. Also, on the boundary, we have

\begin{center}
 $\frac{\partial z_{1}}{\partial \eta}+\sqrt[]{\lambda +\mu} z_{1} =0$ and $\frac{\partial z_{2}}{\partial \eta}+\sqrt[]{\lambda +\mu} z_{2} =0$.
\end{center}
 Hence  ($  z_{1} $, $ z_{2}$) is a supersolution of (1) for $M_{\lambda,\mu}\gg1$.\\\\
$(ii)$ \noindent (F2), (F3) are mirror image cases. So, it is enough to consider one case.\\
\noindent  Assume that (F2) is satisfied:
Let
$$( z_{1} , z_{2}) = \Big(  M_{\lambda ,\mu}\zeta_{\lambda,\mu } ,  \lambda g(M_{\lambda ,\mu}({\|\zeta_{\lambda,\mu }\|}_\infty+M_{\lambda ,\mu}) \zeta_{\lambda,\mu }\Big)$$

\noindent We want to show that
\begin{center}
  $-\Delta z_{2} \geq \lambda g(z_{1})+\mu q(z_{2}).$
\end{center}

\noindent By mean value theorem, there exists $A\in (M_{\lambda ,\mu}{\|\zeta_{\lambda,\mu }\|}_\infty, M_{\lambda ,\mu}{\|\zeta_{\lambda,\mu }\|}_\infty+M_{\lambda ,\mu})$ such that

$$\lambda M_{\lambda ,\mu}g'(A)=\lambda g(M_{\lambda ,\mu}{\|\zeta_{\lambda,\mu }\|}_\infty+M_{\lambda ,\mu} ) - \lambda g(M_{\lambda,\mu }{\|\zeta_{\lambda,\mu }\|}_\infty).$$
Assuming that  $M_{\lambda ,\mu}\gg 1$ (ie. $A\gg 1$), we have
\begin{eqnarray*}
\begin{split}
\mu q(\lambda g(M_{\lambda ,\mu}{\|\zeta_{\lambda,\mu }\|}_\infty+M_{\lambda ,\mu})  \zeta_{\lambda,\mu })&< \mu (\lambda \zeta_{\lambda,\mu })^\beta (g(M_{\lambda ,\mu}{\|\zeta_{\lambda,\mu }\|}_\infty+M_{\lambda ,\mu}))^\beta \\
&< \mu (\lambda \zeta_{\lambda,\mu })^\beta ((M_{\lambda ,\mu}{\|\zeta_{\lambda,\mu }\|}_\infty+M_{\lambda ,\mu})^\gamma))^\beta \\
&=\mu (\lambda \zeta_{\lambda,\mu })^\beta (M_{\lambda ,\mu})^{\gamma \beta}({\|\zeta_{\lambda,\mu }\|}_\infty+1)^{\gamma \beta}<\lambda M_{\lambda ,\mu}g'(A)
\end{split}
\end{eqnarray*}
\noindent for $M_{\lambda ,\mu}\gg1$ since $\gamma \beta<1$ and $\lim\limits_{s\rightarrow\infty}g'(s)>0$.
then we have
$$
\lambda g(M_{\lambda ,\mu}{\|\zeta_{\lambda,\mu }\|}_\infty+M_{\lambda ,\mu} ) - \lambda g(M_{\lambda,\mu }{\|\zeta_{\lambda,\mu }\|}_\infty)=\lambda M_{\lambda ,\mu}g'(A) > \mu q(\lambda g(M_{\lambda ,\mu}{\|\zeta_{\lambda,\mu }\|}_\infty+M_{\lambda ,\mu})  \zeta_{\lambda,\mu }).
$$
This implies that
$$
\lambda g(M_{\lambda ,\mu}{\|\zeta_{\lambda,\mu }\|}_\infty+M_{\lambda ,\mu} )\geq \lambda g(M_{\lambda,\mu }\zeta_{\lambda ,\mu}) + \mu q(\lambda g(M_{\lambda ,\mu}{\|\zeta_{\lambda,\mu }\|}_\infty+M_{\lambda ,\mu})  \zeta_{\lambda,\mu })
$$
Thus we have
\begin{center}
  $-\Delta z_{2} \geq \lambda g(z_{1})+\mu q(z_{2}).$
\end{center}
To prove another inequality, let
$$( z_{1} , z_{2}) = \Big(  M_{\lambda ,\mu}\zeta_{\lambda,\mu } ,  \lambda g(M_{\lambda ,\mu}({\|\zeta_{\lambda,\mu }\|}_\infty+1) \zeta_{\lambda,\mu }\Big)$$
\noindent We want to show that
\begin{center}
  $-\Delta z_{1} \geq \lambda f(z_{2})+\mu h(z_{1}).$
\end{center}

By combine sublinearity $\lim\limits_{s\rightarrow\infty}\frac{f(M(g(s))}{s}=0$ and sublinearity $\lim\limits_{s\rightarrow\infty}\frac{h(s)}{s}=0$ we have

\begin{eqnarray*}
\begin{split}
\frac{1}{{\|\zeta_{\lambda ,\mu}\|}_\infty}&\geq \frac{ \lambda f( \lambda g(M_{\lambda ,\mu}({\|\zeta_{\lambda,\mu }\|}_\infty+1)+\mu h(M_{\lambda ,\mu}\zeta_{\lambda ,\mu})}{ M_{\lambda ,\mu}{\|\zeta_{\lambda ,\mu}\|}_\infty}\\
&=\frac{ \lambda f( \lambda g(M_{\lambda ,\mu}({\|\zeta_{\lambda,\mu }\|}_\infty+1)+\mu h(M_{\lambda ,\mu}\zeta_{\lambda ,\mu})}{ M_{\lambda ,\mu}({\|\zeta_{\lambda,\mu }\|}_\infty+1)}.\frac{M_{\lambda ,\mu}({\|\zeta_{\lambda,\mu }\|}_\infty+1)}{M_{\lambda ,\mu}{\|\zeta_{\lambda ,\mu}\|}_\infty}+\frac{\mu h(M_{\lambda ,\mu}{\|\zeta_{\lambda ,\mu}\|}_\infty)}{M_{\lambda ,\mu}{\|\zeta_{\lambda ,\mu}\|}_\infty}.
\end{split}
\end{eqnarray*}
This implies that
$$
M_{\lambda ,\mu}\geq  \lambda f( \lambda g(M_{\lambda ,\mu}({\|\zeta_{\lambda,\mu }\|}_\infty+1)\zeta_{\lambda ,\mu})+\mu h(M_{\lambda ,\mu}\zeta_{\lambda ,\mu})
$$
Thus we have
\begin{center}
  $-\Delta z_{1} \geq \lambda f(z_{2})+\mu h(z_{1}).$
\end{center}
Further, on the boundary, we have
 \begin{center}
$\frac{\partial z_{1}}{\partial \eta}$+$\sqrt[]{\lambda +\mu} z_{1} $= 0  and
$\frac{\partial z_{2}}{\partial \eta}$+$\sqrt[]{\lambda +\mu}  z_{2}$= 0.
\end{center}
 Hence  ($  z_{1} $, $ z_{2}$) is a supersolution of (1) for  $M_{\lambda,\mu }\gg1$.\\

(iii) Finally we consider the case where $f(x)$$ \rightarrow \infty$ and $g(x)$$ \rightarrow \infty$ as ${x}  \rightarrow \infty$ and $h$ and $q$ are bounded In this case, let
\begin{center}
 $(  z_{1} ,  z_{2})= (\lambda f(M_{\lambda, \mu }{\|\zeta_{\lambda,\mu }\|}_\infty )\zeta_{\lambda ,\mu} ,  \lambda g(M_{\lambda, \mu }{\|\zeta_{\lambda,\mu }\|}_\infty )\zeta_{\lambda ,\mu} ).$
\end{center}
Since $f(x)\rightarrow \infty$ as ${x}  \rightarrow \infty,$ by a Similar argument we see that
 $$
 \lambda f(M_{\lambda ,\mu}{\|\zeta_{\lambda,\mu }\|}_\infty )\geq \lambda f(g(M_{\lambda, \mu }{\|\zeta_{\lambda,\mu }\|}_\infty )\zeta_{\lambda ,\mu}) + \mu h(\lambda f(M_{\lambda,\mu }{\|\zeta_{\lambda ,\mu}\|}_\infty) \zeta_{\lambda,\mu }) .
 $$
 This implies
 \begin{center}
 $-\Delta z_{1} \geq\lambda f(z_{2})+\mu h(z_{1})$.
\end{center}
Similarly, Since $g(x)\rightarrow \infty$ as ${x}  \rightarrow \infty$, we have
 $$
 \lambda g(M_{\lambda ,\mu}{\|\zeta_{\lambda,\mu }\|}_\infty )\geq \lambda g(f(M_{\lambda, \mu }{\|\zeta_{\lambda,\mu }\|}_\infty )\zeta_{\lambda ,\mu}) + \mu q(\lambda g(M_{\lambda,\mu }{\|\zeta_{\lambda ,\mu}\|}_\infty) \zeta_{\lambda,\mu }) .
 $$
 This implies
 \begin{center}
$-\Delta z_{2} \geq \lambda g(z_{1})+\mu q(z_{2}). $
\end{center}
 Further, on the boundary, we have
\begin{center}
 $\frac{\partial z_{1}}{\partial \eta}$+$\sqrt[]{\lambda +\mu} z_{1} $ = 0  and
$\frac{\partial z_{2}}{\partial \eta}$+$\sqrt[]{\lambda +\mu}  z_{2}$ = 0.
\end{center}
Hence  ($  z_{1} $, $ z_{2}$) is a supersolution of (1) for  $M_{\lambda ,\mu}\gg1$. Now  we can choose   $M_{\lambda ,\mu}$  lare enough such that ($  \psi_{1} $, $ \psi_{2}$) $\leq$ ($  z_{1} $, $ z_{2}$). Hence (1) has a positive solution for  ${\lambda+\mu }>A_1$.\\

Next, we show that there exist a positive solution ($ u _{\lambda,\mu} $, $ v_{\lambda ,\mu}$) for ${\lambda+\mu }\gg1$
 such that ${\|  u _{\lambda,\mu }  \|}$, ${\|  v _{\lambda ,\mu}  \|} \rightarrow \infty$ as ${\lambda +\mu}  \rightarrow \infty$.

Define $k{\in C^2 [0,\infty)}$  such that $k(0) < 0$, $k(s)\leq {f(s)}, { h(s)}, {g(s)}$ and  ${q(s)}$ for $s\in(0,\infty)$ and   $\lim\limits_{s\to \infty}{k(s) > 0}  $. Then the Dirichlet boundary value problem

\begin{equation}
\left\{\begin{array}{ll}
-\Delta P= (\lambda+\mu) k(P);  &   \Omega,\\
P=0; & \partial\Omega,    \\
\end{array}\right.
\end{equation}
 has a positive solution  $\overline {P} _{\lambda,\mu }$   for  ${\lambda +\mu}\gg1$ such that  ${\|  \overline {P} _{\lambda ,\mu}  \|} \rightarrow \infty$ as  ${\lambda +\mu}  \rightarrow \infty$ (see \cite{ac-jbg-rs}).

 It is easy to show that ($\overline {P} _{\lambda,\mu }$ , $\overline {P} _{\lambda ,\mu}$ ) is a subsolution of (1) for ${\lambda +\mu}\gg1$ since $k\leq f+ h$, $k\leq g+\gamma$, and
$\frac{\partial \overline {P} _{\lambda ,\mu}}{\partial \eta} <0 ; \partial\Omega$ by Hoop maximum principle. We can also choose  $M_{\lambda,\mu }\gg1$ such that   $(  z_{1} , z_{2}) \geq(\overline {P} _{\lambda ,\mu},\overline {P} _{\lambda,\mu } )$,  where   ($  z_{1} $, $ z_{2}$) is a supersolution of (1) as constructed before.

This implies that (1) has a positive solution  ($u _{\lambda,\mu }$, $v _{\lambda ,\mu}$)   such that

\begin{center}
  $(u _{\lambda,\mu } ,v _{\lambda,\mu })\in[ (\overline {P} _{\lambda ,\mu},\overline {P} _{\lambda ,\mu}) ,  (  z_{1} ,z_{2}) ]$ for ${\lambda +\mu}\gg1$.
\end{center}
Clearly, ${\|  u _{\lambda ,\mu}  \|}$, ${\|  v _{\lambda ,\mu}  \|} \rightarrow \infty$ as ${\lambda+\mu }  \rightarrow \infty$  since    ${\|  \overline {P} _{\lambda,\mu }  \|} \rightarrow \infty$ as  ${\lambda+\mu }  \rightarrow \infty$.

Next, we establish our multiplicity result. We first construct a positive strict supersolution ($ \tilde{ z_{1}} $, $\tilde{ z_{2}}$) for (1) when $\lambda ,\mu \in (1, \frac{Q_1}{2\|\xi\|}_\infty)$, where $\xi$ is the solution of

\begin{equation}
\left\{\begin{array}{ll}
-\Delta \xi =1;  &   \Omega,\\
\frac{\partial \xi}{\partial \eta}+ \xi =0; & \partial\Omega,    \\
\end{array}\right.
\end{equation}
Let  ($ \tilde{ z_{1}} $, $\tilde{ z_{2}}$) =($a \frac{\xi}{\|\xi\|}_\infty$, $a \frac{\xi}{\|\xi\|}_\infty$).
We have

\begin{center}
$\lambda <\frac{Q_1 }{2\|\xi\|}_\infty\leq\frac{1}{2\|\xi\|}_\infty \frac{a}{f(a)}.$
\end{center}
Then
\begin{equation}
\frac{a}{2\|\xi\|}_\infty >\lambda f(a)
\geq \lambda f\Big(a\frac{\xi}{\|\xi\|}_\infty \Big),
\end{equation}

(since $f$ is increasing). Also since  $\mu <\frac{Q_1}{2\|\xi\|}_\infty,$ similarly we obtain

\begin{center}
$\mu <\frac{Q_1 }{2\|\xi\|}_\infty\leq\frac{1}{2\|\xi\|}_\infty \frac{a}{h(a)}.$
\end{center}
Then
\begin{equation}
\frac{a}{2\|\xi\|}_\infty >\mu h(a)
\geq \mu h\Big(a\frac{\xi}{\|\xi\|}_\infty\Big).
\end{equation}
From $(18)$ and $(19)$ we have $$\frac{a}{\|\xi\|}_\infty>\lambda f\Big(a\frac{\xi}{\|\xi\|}_\infty \Big)+\mu h\Big(a\frac{\xi}{\|\xi\|}_\infty \Big),$$ which implies that $-\Delta \tilde z_{1} \geq \lambda f(\tilde z_{2})+\mu h(\tilde z_{1}).$

Similarly, we can show that $-\Delta \tilde z_{2} \geq \lambda g(\tilde z_{1})+\mu q(\tilde z_{2})$  in $\Omega$. On the boundary, we have
\begin{eqnarray*}
\begin{split}
\frac{\partial \tilde z_{1}}{\partial \eta}+\sqrt[]{\lambda +\mu} \tilde z_{1} &= \frac{a}{\|\xi\|}_\infty\frac{\partial \xi}{\partial \eta}+\sqrt[]{\lambda +\mu} \frac{a}{\|\xi\|}_\infty \xi  \\
&=\frac{a}{\|\xi\|}_\infty(\frac{\partial \xi}{\partial \eta}+\sqrt[]{\lambda +\mu} \xi) \\
&> \frac{a}{\|\xi\|}_\infty(\frac{\partial \xi}{\partial \eta}+ \xi)\\
&= 0.
\end{split}
\end{eqnarray*}
This implies $\frac{\partial \tilde z_{1}}{\partial \eta}+\sqrt[]{\lambda +\mu} \tilde z_{1} > 0$.

Similarly, we can show that $\frac{\partial \tilde z_{2}}{\partial \eta}+\sqrt[]{\lambda +\mu} \tilde z_{2} >0$. Thus  ($ \tilde{ z_{1}} $, $\tilde{ z_{2}}$) is a strict supersolution of (1).

Now we construct a strict subsolution  ($ \tilde \psi_{1} $,$\tilde \psi_{2}$)  of (1) for $\lambda , \mu\geq C_1Q_2$.
Note that in \cite{ja-mr-rs} the authors showed that the system

\begin{equation}
\left\{\begin{array}{ll}
-\Delta u = \lambda f(v);  &   \Omega,\\
-\Delta  v=  \lambda g(u);& \Omega,\\
u=0; & \partial\Omega,\\
v=0; & \partial\Omega,\\
\end{array}\right.
\end{equation}
has a strict subsolution $(\overline{u_0} , \overline{v_0})$ for  $\lambda \geq C_1 $max$\{\frac{b}{f(b)}, \frac{b}{g(b)}\}$ such that ${\|\overline{u_0} \|}_\infty \geq b$ and ${\|\overline{v_0} \|}_\infty \geq b$. Let   ($ \tilde \psi_{1} $, $\tilde \psi_{2}$)   be the first iteration of $(\overline {u_0} , \overline {v_0})$, that is, the solution to the problem

\begin{equation}
\left\{\begin{array}{ll}
-\Delta \tilde \psi_1 = \lambda f(\overline{v_0})+\mu h(\overline{u_0}) ;  &   \Omega,\\
-\Delta    \tilde \psi_2=  \lambda g(\overline{u_0})+\mu q(\overline{v_0}) ;& \Omega,\\
\frac{\partial \tilde \psi_1}{\partial \eta}+\sqrt[]{\lambda +\mu} \tilde \psi_1=0; & \partial\Omega,\\
\frac{\partial \tilde \psi_2}{\partial \eta}+\sqrt[]{\lambda +\mu} \tilde \psi_2=0; & \partial\Omega.\\
\end{array}\right.
\end{equation}
Then by comparison principle  ($ \tilde \psi_{1} $, $\tilde \psi_{2}$) $ >$ $(\overline {u_0} , \overline {v_0})$; $\Omega$. Hence ($ \tilde \psi_{1} $, $\tilde \psi_{2}$) is a strict subsolution of (1) such that ${\|\tilde \psi_1 \|}_\infty \geq b>a$  and  ${\|\tilde \psi_2 \|}_\infty \geq b>a$. Thus we have ($ \tilde \psi_{1} $, $\tilde \psi_{2}$)$\nleq$  ($ \tilde{ z_{1}} $, $\tilde{ z_{2}}$). Recall the subsolution ($ \psi_{1} $, $ \psi_{2}$)  for  $ m\approx 0$ when $\lambda+\mu >A_{1}$ and the supersolution  ($ z_{1} $, $ z_{2}$) for $M_{\lambda ,\mu}\gg1$ when $\lambda, \mu>0$. Also, for $ m\approx o $ and  $M_{\lambda ,\mu}\gg1$, we have $(\psi_{1}, \psi_{2}) \leq (\tilde \psi_{1} , \tilde \psi_{2}) \leq ( z_{1} ,  z_{2})$ and $(\psi_{1}, \psi_{2}) \leq  (\tilde z_{1} ,\tilde z_{2}) \leq ( z_{1} ,  z_{2}).$ Hence by Lemma $(3.1)$, $(1)$ has at least three positive solution for $\lambda ,\mu \in (\max \{A_1, C_1Q_2,1\},\frac{Q_1}{2\|\xi\|}_\infty)$. This completes the proof.\\\\
{\bf Proof of Theorem 2.2}. Note that $(\psi_{1},\psi_{2}) = ( m{\phi} _{\lambda ,  \mu} , m {\phi} _{\lambda ,  \mu} )$ is a subsolution of $(1)$ for  $ m\approx 0.$ Since, $  f^{\prime\prime}(s), g^{\prime\prime}(s), h^{\prime\prime}(s)$ and $ q^{\prime\prime}(s) < 0$  for  $ s\approx 0$, there exists  $ M > 0 $ such that  $f^{\prime\prime}(s),  g^{\prime\prime}(s),  h^{\prime\prime}(s)$ and $q^{\prime\prime}(s)\leq-M$ for  $s\approx 0$.

Let ($ \phi_{1} $, $ \phi_{2}$) $ = $ $( \vartheta _{\lambda ,  \mu}{\phi} _{\lambda ,  \mu},\vartheta _{\lambda ,  \mu} {\phi} _{\lambda ,  \mu} )$, where $ \vartheta _{\lambda ,  \mu} = -\frac {2 \varrho _{\lambda ,  \mu}}{ (\lambda +  \mu)M{\min}_{\overline{\Omega}}  \phi_{\lambda ,  \mu}}$. Note that $ \vartheta _{\lambda ,  \mu}>0$ and $ \vartheta _{\lambda ,  \mu}  \rightarrow 0$ as ${\lambda +  \mu}  \rightarrow A_1^+$ and  $ {{\min}_{\overline{\Omega}}  \phi_{\lambda ,  \mu}}  \nrightarrow 0$  as ${\lambda +  \mu}  \rightarrow A_1^+$ (see \cite{jg-qm-sr-rs}). This implies that ${\| \phi_1\|}_\infty \rightarrow 0$ and ${\| \phi_2\|}_\infty   \rightarrow 0$ as ${\lambda +  \mu}  \rightarrow A_1^+$.
Then by Taylors theorem we have
$$
f(\phi_2)= f(0) + f^\prime(0)\phi_2 + \frac{ f^{\prime\prime}(\zeta)}{2}\phi_2 ^2 ,
$$
for some $\zeta\in[0,\phi_2 ]$, and
$$
h(\phi_1)= h(0) + h^\prime(0)\phi_1 + \frac{ h^{\prime\prime}(\zeta)}{2}\phi_1 ^2 ,
$$
for some $\zeta\in[0,\phi_1].$ Then we have

\begin{eqnarray*}
\begin{split}
&-\Delta  \phi_1 -\lambda f(\phi_2) -\mu h(\phi_1)&\\
&=  \vartheta _{\lambda ,  \mu}(\varrho_{\lambda ,  \mu} +(\lambda +\mu) g^\prime(0)) \phi_{\lambda ,\mu}
-\lambda [ \vartheta _{\lambda ,  \mu} \phi_{\lambda ,  \mu}f^\prime(0) + \frac{ f^{\prime\prime}(\zeta)}{2}{  (\vartheta _{\lambda ,  \mu} \phi_{\lambda ,  \mu})} ^2]\\
&-\mu [ \vartheta _{\lambda ,  \mu} \phi_{\lambda ,  \mu}h^\prime(0) + \frac{ h^{\prime\prime}(\zeta)}{2}{  (\vartheta _{\lambda ,  \mu} \phi_{\lambda ,  \mu})} ^2] \\
&=\vartheta _{\lambda ,  \mu}  \varrho_{\lambda ,  \mu} \phi_{\lambda ,  \mu} +\vartheta _{\lambda ,  \mu}  \phi_{\lambda ,  \mu}(\lambda +\mu) g^\prime(0)
-\lambda \vartheta _{\lambda ,  \mu}  \phi_{\lambda ,  \mu}f^\prime(0)-  \frac{\lambda f^{\prime\prime}(\zeta)}{2}{
(\vartheta _{\lambda ,  \mu} \phi_{\lambda ,  \mu})} ^2 \\
&-\mu \vartheta _{\lambda ,  \mu}  \phi_{\lambda ,  \mu}h^\prime(0)-  \frac{\mu h^{\prime\prime}(\zeta)}{2}{  (\vartheta _{\lambda ,  \mu} \phi_{\lambda ,  \mu})} ^2 \\
&=\vartheta _{\lambda ,  \mu}  \varrho_{\lambda ,  \mu} \phi_{\lambda ,  \mu}  +\vartheta _{\lambda ,  \mu}  \phi_{\lambda ,  \mu}(\lambda +\mu) g^\prime(0)
 - \vartheta _{\lambda ,  \mu}  \phi_{\lambda ,  \mu}(\lambda +\mu) g^\prime(0)\\
&-  \frac{\lambda f^{\prime\prime}(\zeta)}{2}{  (\vartheta _{\lambda ,  \mu} \phi_{\lambda ,  \mu})} ^2 -  \frac{\mu h^{\prime\prime}(\zeta)}{2}{  (\vartheta _{\lambda ,  \mu} \phi_{\lambda ,  \mu})} ^2 \\
&\geq \vartheta _{\lambda ,  \mu} \varrho_{\lambda ,  \mu}  \phi_{\lambda ,  \mu}+\frac{\lambda M}{ 2 }{  (\vartheta _{\lambda ,  \mu} \phi_{\lambda ,  \mu})} ^2
+ \frac{\mu M}{ 2 }{  (\vartheta _{\lambda ,  \mu} \phi_{\lambda ,  \mu})} ^2 \\
&= \vartheta _{\lambda ,  \mu} \varrho_{\lambda ,  \mu}  \phi_{\lambda ,  \mu}+\frac{(\lambda + \mu)M}{ 2 }{  (\vartheta _{\lambda ,  \mu} \phi_{\lambda ,  \mu})} ^2 \\
&= \vartheta _{\lambda ,  \mu} \phi_{\lambda ,  \mu}  [\varrho_{\lambda ,  \mu} +\frac{(\lambda + \mu)M}{ 2 }  \vartheta _{\lambda ,  \mu} \phi_{\lambda ,  \mu}] \\
&\geq \vartheta _{\lambda ,  \mu} \phi_{\lambda ,  \mu}  [\varrho_{\lambda ,  \mu} +\frac{(\lambda + \mu)M}{ 2 }  \vartheta _{\lambda ,  \mu} {\min}_{\overline{\Omega}}  \phi_{\lambda ,  \mu} ] \\
&= \vartheta _{\lambda ,  \mu} \phi_{\lambda ,  \mu}  [\varrho_{\lambda ,  \mu} -\frac{(\lambda + \mu)M {\min}_{\overline{\Omega}}  \phi_{\lambda ,  \mu} }{ 2 } \times \frac {2 \varrho _{\lambda ,  \mu}}{ (\lambda +  \mu)M {\min}_{\overline{\Omega}}  \phi_{\lambda ,  \mu}}] \\
&= \vartheta _{\lambda ,  \mu} \phi_{\lambda ,  \mu}  [\varrho_{\lambda ,  \mu} - \varrho_{\lambda ,  \mu} ] \\
&= 0.
\end{split}
\end{eqnarray*}

Then we have

\begin{center}
$-\Delta  \phi_1 -\lambda f(\phi_2) -\mu h(\phi_1) \geq 0,$
\end{center}
by our choose of $\vartheta _{\lambda ,  \mu}$. A similar argument will show that
\begin{center}
 $-\Delta  \phi_2 -\lambda g(\phi_1) -\mu q(\phi_2)\geq 0$
\end{center}
Further, on the boundary we have

\begin{center}
$\frac{\partial  \phi_1}{\partial \eta}+\sqrt[]{\lambda +\mu}  \phi_1=0 $ and $ \frac{\partial  \phi_2}{\partial \eta}+\sqrt[]{\lambda +\mu}  \phi_2=0.$
\end{center}
Thus  ($ \phi_{1} $, $ \phi_{2}$) is a supersolution of $(1).$ Choosing $m\approx 0$, we also have ($ \psi_{1} $, $ \psi_{2}$)$\leq$($ \phi_{1} $, $ \phi_{2}$).

 Therefore, there exists a positive solution $(u _{\lambda , \mu} ,v _{\lambda ,  \mu})\in[ ( {\psi} _{1}, {\psi} _{2}) ,  (  \phi_{1} ,\phi_{2}) ]$ for ${\lambda +  \mu}> A_1$
 such that ${\|  u _{\lambda ,  \mu}  \|}$, ${\|  v _{\lambda ,  \mu}  \|} \rightarrow 0$ as ${\lambda +  \mu}  \rightarrow {A_1}^+$. This completes the proof.

\section{Example}
\hspace{0.6 cm}For an example to illustrate Theorems $2.1$ and $2.2$, consider $f=f_{a,k},h=h_{a , k}, g=g_{k},$ and $q=q_{k}$ as follows:

\begin{equation}
f=f_{a , k}(s)=\begin{cases}
e^\frac{s}{s+1} -1 ;  &   s\leq k,\\
(e^\frac{\alpha s}{ \alpha+ s} -e^\frac{\alpha k}{ \alpha+ k} )  + (e^\frac{k}{k+1} -1) ;& s\geq k,\\
\end{cases}
\end{equation}

\begin{equation}
h=h_{a , k}(s)=\begin{cases}
e^\frac{2s}{s+1}- s -1 ;  &   s\leq k,\\
(e^\frac{\alpha s}{ \alpha+ s} -e^\frac{\alpha k}{ \alpha+ k} )  + (e^\frac{k}{k+1} -1) ;& s\geq k,\\
\end{cases}
\end{equation}

\begin{equation}
g=g_{k}(s)=\begin{cases}
2(1+ s)^\frac{1}{2} -2 ;  &   s\leq k,\\
\Big(\frac{1}{ 2} (1+ s)^2 -\frac{1}{ 2} (1+ k)^2  \Big)  + \Big(2(1+ k)^\frac{1}{2} -2\Big) ;& s\geq k,\\
\end{cases}
\end{equation}
and
\begin{equation}
q=q_{k}(s)=\begin{cases}
3(1+ s)^\frac{1}{3} -3 ;  &   s\leq k,\\
\Big(\frac{1}{ 3} (1+ s)^3 -\frac{1}{ 3} (1+ k)^3  \Big)  + \Big(3(1+ k)^\frac{1}{3} -3\Big) ;& s\geq k,\\
\end{cases}
\end{equation}
where $k>0$ is a constant, and $\alpha >0$ is a parameter. Note that though $g$ and $q$ are superlinear at infinity, since $f$ is bounded, $\frac{f(M(g(s))}{s}\rightarrow 0$ as $s\rightarrow\infty$ for all $M>0$. Choose $a=k$, $b=\alpha$ and $\alpha >k$. Then we note that $Q_{1}(k)=\min\{\frac{k}{f(k)}, \frac{k}{g(k)},\frac{k}{h(k)},\frac{k}{q(k)}\} \rightarrow\infty$ as $k\rightarrow\infty$ since $\frac{k}{f(k)}, \frac{k}{g(k)},\frac{k}{h(k)},\frac{k}{q(k)} \rightarrow\infty$ as $k\rightarrow\infty$ and $Q_{2}(k)=\max\{\frac{\alpha}{f(\alpha)}, \frac{\alpha}{g(\alpha)},\frac{\alpha}{h(\alpha)},\frac{\alpha}{q(\alpha)}\} \rightarrow 0$ as $\alpha\rightarrow\infty$ since $\frac{\alpha}{f(\alpha)}, \frac{\alpha}{g(\alpha)},\frac{\alpha}{h(\alpha)},\frac{\alpha}{q(\alpha)} \rightarrow 0$ as $\alpha\rightarrow\infty$. Hence we can choose $k=k_{0}$ such that $Q_{1}(k_{0})>\max 2\{A_1,1\}{\|\xi\|}_\infty$ and choose $\alpha \gg1$ such that $\frac{Q_{1}(k_{0})}{Q_{2}(\alpha)}>2C_1{\|\xi\|}_\infty$. Hence  $\frac{Q_{1}(k_{0})}{2{\|\xi\|}_\infty} >\max \{A_1, C_{1}Q_{2},1\}$ for $k=k_{0}$ and $\alpha \gg1$. It is also clear that $f, g, h$ and $q$ satisfy $(H_{1})- (H_{3})$. This implies that $(1)$ has at least one positive solution for $\lambda + \mu >A_{1}$ and at least three positive solution for $\lambda ,\mu \in (\max \{A_1, C_{1}Q_{2},1\},\frac{Q_{1}}{{2\|\xi\|}_\infty})$ for $k=k_{0}$ and $\alpha \gg1$. Thus Theorem $2.1$ holds in this example for $k=k_{0}$ and $\alpha \gg1$. Further, since $f^\prime(0) = g^\prime(0) = h^\prime(0) = q^\prime(0) > 0$ and $f, g \in C^{2} [0,k)$  and $ f^{\prime\prime}(s) , g^{\prime\prime}(s), h^{\prime\prime}(s) , q^{\prime\prime}(s) < 0$ for all $ s\in [0, k)$, Theorem $2.2$ holds.\\\\

\section*{Declarations}

\subsection*{Ethical Approval }

Hereby, We consciously assure that the following is fulfilled:\\

$1)$ This material is the authors' own original work, which has not been previously published elsewhere.\\

$2)$ The paper is not currently being considered for publication elsewhere.\\

$3)$ The paper reflects the authors' own research and analysis in a truthful and complete manner.

\subsection*{Competing interests}

The authors state no conflict of interest.
\subsection*{ Availability of data and materials}
Data sharing not applicable to this article as no datasets were generated or analyzed during the current study.
\subsection*{Funding statement}
There are no funders to report for this submission.

\end{document}